\documentclass[a4paper,10pt]{article}
\usepackage{amsmath,amsthm}
\usepackage{amssymb}
\usepackage{psfrag}
\usepackage[dvips]{graphicx}
\usepackage{placeins}

\newtheorem{theorem}{Theorem}
\newtheorem{proposition}{Proposition}
\newtheorem{lemma}{Lemma}

\def\RR{\mathbb{R}}
\def\ZZ{\mathbb{Z}}

\def\dim{\mathrm{dim}}
\def\WB{\mathit{WB}}

\psfrag{ac}[][b1][0.6]{$a\sqrt{2}$}
\psfrag{bc}[][b1][0.6]{$b\sqrt{2}$}
\psfrag{c.}[][b1][0.6]{$\sqrt{2}$}
\psfrag{p1}[][b1][0.8]{$p_1$}
\psfrag{p2}[][b1][0.8]{$p_2$}
\psfrag{p3}[][b1][0.8]{$p_3$}
\psfrag{p4}[][b1][0.8]{$p_4$}

\title{A special tiling of the rectangle}
\author{Carlos Tomei and Tania Vieira \footnote{Departamento de Matem\'atica, PUC-Rio;
tomei@mat.puc-rio.br, tvieira@mat.puc-rio.br.
}}
\date{}
\begin{document}
\maketitle

\begin{abstract}
We count tilings of a rectangle of integer sides $m-1$ and $n-1$
by a special set of tiles. The result is obtained from the study of the kernel of the adjacency
matrix of an $m \times n$ rectangular subgraph in $\ZZ \times \ZZ$.
\end{abstract}

\section{Introduction}

\begin{figure}[hbtp]
\centering
\includegraphics[height=1.2in]{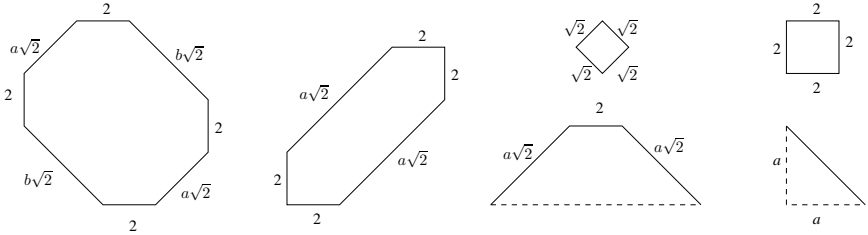}
\caption{The tiles; $a$ and $b$ are integers.}
\label{pecas}
\end{figure}

\begin{figure}[hbtp]
$$
\begin{array}{ccc}
\includegraphics[width=1.3in]{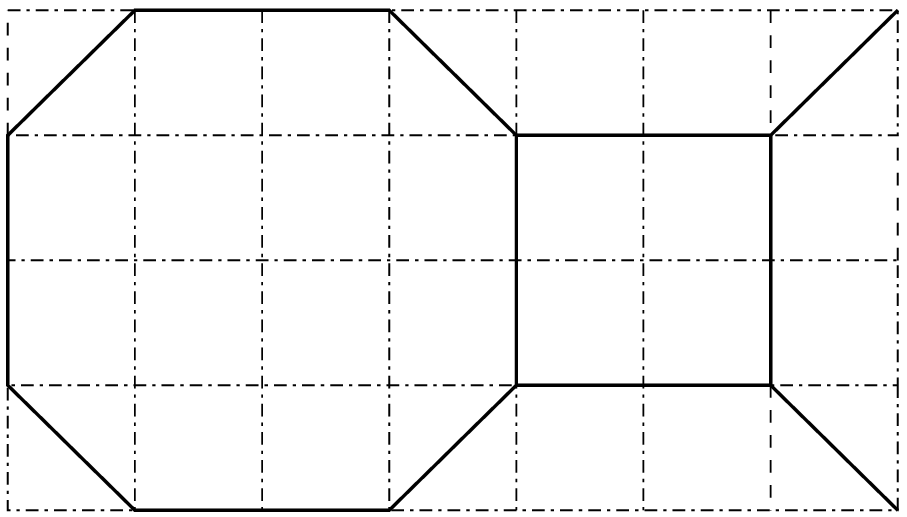}&\hspace*{0.5cm}&
\includegraphics[width=1.3in]{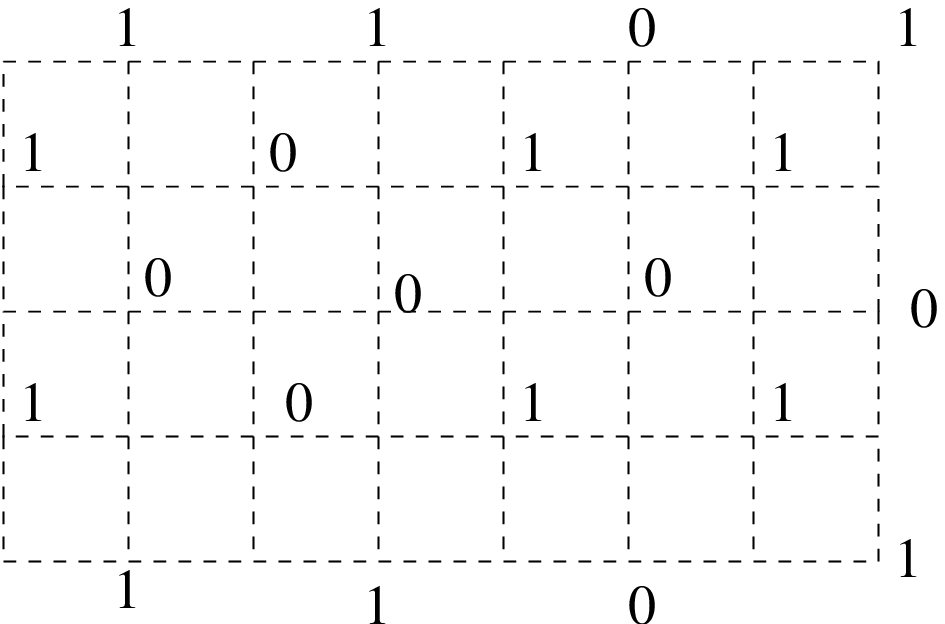}\\\\
( 1 )&&( 2 )
\end{array}
$$
\caption{A tiling and a polarized $\ZZ/(2)$-harmonic function of the $7 \times 4$ rectangle.}
\label{exemplo de ladrilhamento} 
\end{figure}

Let $R = [0,m-1] \times [0,n-1] \subset\RR ^2$ be a rectangle with nontrivial integer sides 
$m-1$ and $n-1$.
We consider tilings of $R$ by tiles given in figure \ref{pecas}. Tiles
can be rotated by integer multiples of $\pi/2$ and their dashed sides ought
to belong to $\partial R$, the boundary of $R$. Length of sides is 
indicated in the picture, angles must be  $\pi/4, \pi/2$ or $3\pi/4$. 
Figure \ref{exemplo de ladrilhamento}.1 shows a tiling of the $7 \times 4$ rectangle.

The main result of this paper is the following.

\begin{theorem} \label{formula}
Let $R = [ 0,m-1 ] \times [ 0,n-1 ]$ be a rectangle with nonzero integer sides $m-1$ and \mbox{$n-1$}.
The number of tilings of $R$ which make use of the tiles given in figure
\ref{pecas} is \mbox{$2^{\beta}+2^{\omega}-2$}, where 
\mbox{$c=\gcd(m+1,n+1)-1$},
\mbox{$2\beta = c+ (c\bmod 2)$} and \mbox{$2\omega = c- (c\bmod 2)$}.
\end{theorem}

The theorem is a consequence of a certain description of 
{\emph{$\ZZ/(2)$-harmonic functions}}. More precisely, consider an $m \times n$ rectangular
mesh of points in $\ZZ \times \ZZ$ and color the points black and white alternatively.
The adjacency matrix $M$ for the points of the mesh is
a standard discretization of the Laplacian acting on functions defined on the rectangle
$R$, satisfying Dirichlet boundary conditions, up to a multiplicative factor given by the 
mesh spacing and to an additive factor given by four times the identity matrix.

The matrix $M$ takes vectors supported on the set of black points (forming a subspace $V_b$)
to vectors supported on white points (in $V_w$), and vice-versa.
Thus, $M$ decomposes naturally in two linear maps, \mbox{$BW: V_b \to V_w$} and
\mbox{$\WB: V_w \to V_b$} and the kernel of $M$ is the direct sum of the 
kernels of $BW$ and $\WB$. 
We assume the field of scalars for both vector spaces to be $\ZZ/(2)$, 
the field of two elements: in this case,
the adjacency matrix equals the discretized Laplacian. The elements on
$\ker BW$ and $\ker \WB$ are the {\emph{polarized}} $\ZZ/(2)$-harmonic functions, 
supported respectively on black and white points. Such objects have already been 
considered in  \cite{Tomei} as a technical tool.

In order to prove the Theorem, we first construct a bijection between tilings of $R$
and the nontrivial (i.e., nonzero) polarized $\ZZ/(2)$-harmonic functions: this 
takes most of Sections 3 and 4. We are then left with counting the dimensions of the 
kernels of $BW$ and $\WB$ for rectangles. 
This in turn is accomplished by realizing that polarized harmonic functions
on $R$ admit very stringent symmetries: the problem then reduces to computing the kernel
dimensions for squares. We conclude by listing the tilings of the $10 \times 4$ rectangle, 
indicating the visual implications of the additive structure on tilings inherited from
the domain vector spaces. 

The first author acknowledges support from CNPq and FAPERJ, the second acknowledges
support from CNPq. Both are thankful to Nicolau Saldanha for many conversations.

\section{Basic properties of tilings}
{\emph{Tiles}} are the polygons (topologically, closed disks) listed in figure
\ref{pecas}. A {\emph{tiling}} of $R = [0,m-1] \times [0,n-1]$ is a covering of $R$ 
by tiles which overlap at their boundaries. Points in the {\emph{$m \times n$ rectangular graph}}
$G=R\cap (\ZZ \times \ZZ)$ are adjacent in the obvious manner: 
$(x,y)$ and $(x^{\prime},y^{\prime})$ are adjacent if and only if
$x=x^{\prime}$, $y=y^{\prime}\pm 1$ or $x=x^{\prime}\pm 1$, $y=y^{\prime}$. 
Notice that we call {\emph{points}} the elements of $G$: {\emph{vertices}} will be the corners of tiles.
Color a point $(x,y) \in G$ {\emph{black}} (resp. {\emph{white}}) if $x+y$ is even (resp. odd).
Points of $\ZZ \times \ZZ$ will be called {\emph{integral points}}.

\begin{proposition}
\label{coordenadas inteiras}
Vertices of tiles of a tiling $T$ are always integral. Also, with the possible exception of 
the four corners of $R$, all vertices of tiles in $T$ are of the same color.
\end{proposition}

\begin{proof}
From the shape of the tiles, if one vertex of a tile sits on an integral point, all
other vertices do. Suppose now that $T$ contains a tile with some nonintegral vertex.
Let $N$ be the region of $R$ covered by tiles of $T$ having nonintegral vertices: $N$ is
a polygonal region which does not contain the corners of $R$ --- corners belong to tiles
with integral vertices. Let $v$ be a vertex of $N$: then it must also be the vertex of some
tile outside $N$ --- a contradiction. Thus, all vertices are integral. 

Again from the shape of tiles, all vertices of a tile have the same color (unless the
tile is a right triangle with small sides of odd length) --- call this the color of the
tile (for triangles, the color of the tile is the color of the vertices not sitting at
the right angle). The right angle vertex of a triangle necessarily sits on a corner of $R$, and we
do not have to consider it in the sequel.

Again, with the possible exception of sides of triangles, all other horizontal or vertical
sides are of even length. Thus, all tiles with a side lying on $\partial R$ have
the same color $c$. Let $C$ be the region covered by tiles of color $c$, and $D$ be its
complement in $R$. If $D$ is not empty, take $v$ to be one of its vertices: since $v \in D$
is not a corner of $R$, it must have the colors of tiles both in $C$ and $D$ --- again a
contradiction.
\end{proof}

Thus $\mathcal{T} = \mathcal{T}_b \cup \mathcal{T}_w$, where $\mathcal{T}$ is the 
set of tilings of $R$, and $\mathcal{T}_b$ 
and $\mathcal{T}_w$ are the
sets of tilings whose tiles are black and white. 

Let $G$ be the rectangular graph associated to $R$. Number
black and white points from $1$ to $n_b$ and from $1$ to $n_w$, respectively.
Set $BW$ to be the $n_w \times n_b$ {\emph{black-to-white adjacency matrix}}, 
with obvious entries: $BW_{wb}=1$ if point 
$b$ is adjacent to point $w$; otherwise, $BW_{wb}=0$. 
Similarly, define $\WB$, the $n_b \times n_w$
{\emph{white-to-black adjacency matrix }} of $G$. 

Figure \ref{exemplo de ladrilhamento}.2 describes  an element in $\ker \WB$ for 
the $7 \times 4$ rectangle (and hence, for the $8 \times 5$ rectangular graph). 
Notice that the nonzero coordinates of this vector correspond to the points of $G$
which belong to the non-dashed sides of the tiling in figure
\ref{exemplo de ladrilhamento}.1: kernel elements naturally induce tilings. This will
be made precise in the next section: we are aiming at the following
result.

\begin{theorem} \label{principal}
There are bijections $$\Phi_b: \ker BW \setminus \{0\} \to \mathcal{T}_b 
\mbox{ and } \Phi_w: \ker \WB \setminus \{0\} \to \mathcal{T}_w$$
\end{theorem}
 
\section{Defining the bijections $\Phi_b$ and $\Phi_w$}

Let $u\in \ker BW \setminus\{0\}$. A black point 
$b \in G$ is an {\emph{active point of $u$}} if $u(b)=1$.
Denote by $A_u$ the set of active points of $u$. 
In $A_u$, consider the following {\emph{adjacency relation}}.
Two points are adjacent if both conditions  are satisfied:
(i) both points have a common white neighbor $w$ in $G$;
(ii) both points form a right angle with $w$ or they are the only 
active neighbors of $w$. 

\begin{figure}[hbtp]
\centering
\includegraphics[height=0.5in]{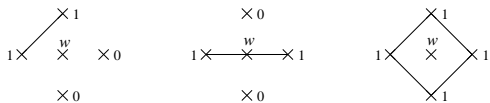}
\caption{Adjacencies in $A_u$.}
\label{adjacencia}      
\end{figure}

A {\emph{splitting}} of $R$ is a decomposition in
{\emph{chunks}}, i.e., polygons with integral vertices, which are not necessarily tiles.
Given a vector $u \in \ker BW \setminus\{0\}$, let 
$\Phi_b(u) = \mathrm{int}(R) \setminus \{ \overline{b_1b_2} \;| 
\mbox{ $b_1$ and $b_2$ are adjacent in $A_u$}\}$.
Said differently, draw straight line segments joining adjacent points in $A_u$.
In section 3.1, we will see that this construction obtains a splitting $\Phi_b(u)$ of $R$. 
In Section 3.3, more will be proved: $\Phi_b(u)$ is actually a tiling, i.e., its chunks are tiles.

\subsection{Geometry of $A_u$}

Let $b \in A_u$ be a black point in the interior of $R$. The 
{\emph {star}} centered at $b$ is the set of points in $A_u$ which are adjacent to $b$, 
together with the segments joining $b$ to these points. The Lemma below follows
by exhausting possibilities.

\begin{figure}[hbtp]
\centering
\includegraphics[height=2.6in]{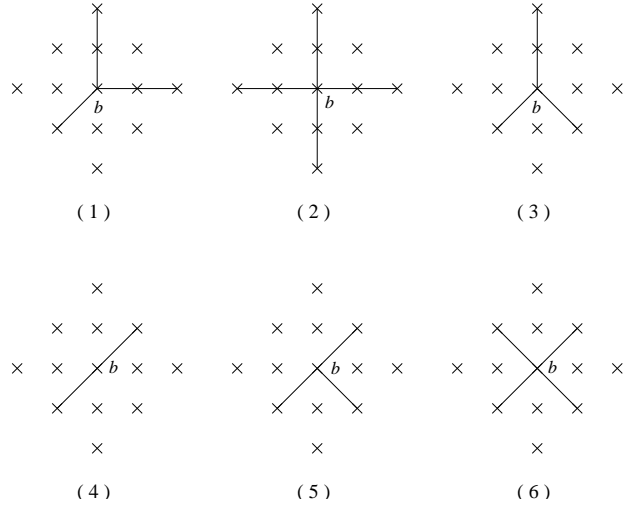}
\caption{Stars centered at $b$.}
\label{estrelas}
\end{figure}

\begin{lemma} 
\label{estrela}
Let $b \in A_u$ be an interior point of $G$. Up to rotations and reflections,
the stars centered at $b$ are listed in figure \ref{estrelas}.
\end{lemma}

We point out that all stars in the list indeed occur in tilings.

\begin{lemma}\label{degree}
Every point of $A_u$ which is not a corner of $R$ has (at least) two neighbors in $A_u$.  
Also, segments $\overline{b_1 b_2}$ and $\overline{b_3 b_4}$ in $A_u$
joining four distinct points do not intercept.

\end{lemma}

\begin{figure}[hbtp]
$$
\begin{array}{ccccc}
\includegraphics[width=1.6in]{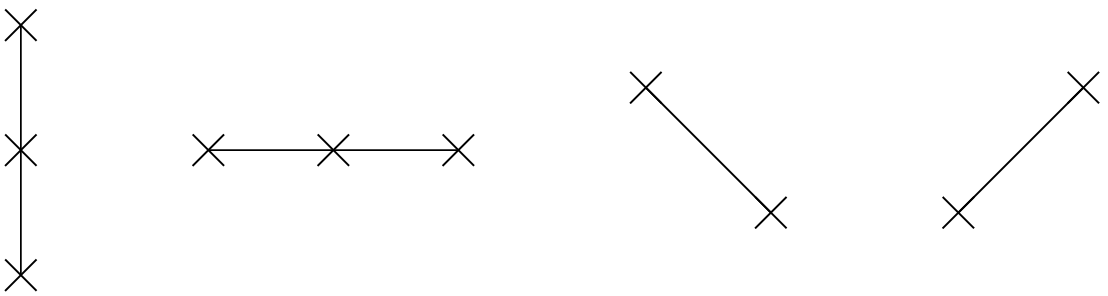}&\hspace*{0.5cm}&
\includegraphics[width=0.45in]{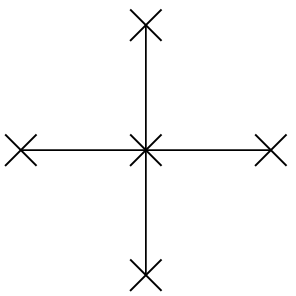}&\hspace*{0.5cm}&
\includegraphics[width=0.25in]{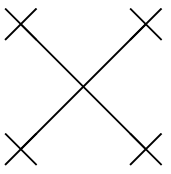}\\\\
\end{array}
$$
\caption{Segments and their possible intersections.}
\label{intersecoes} 
\end{figure}

\begin{proof}
Let $b$ be a point in $A_u$ which is not a corner of $R$.
Then $b$ has two neighbors in $G$, $w_1$ and $w_2$, so that
the three points are collinear. Since $b$ is an active neighbor of 
$w_1$ and $w_2$ and $u\in \ker BW$, 
$w_1$ and $w_2$ have other active black neighbors in $G$. 
Let $b_1$ (resp. $b_2$) be the active neighbor of $w_1$ (resp. $w_2$) closest to $b$.
Notice that $b_1$ and $b_2$ are distinct, since
$b$ is the only black neighbor common to $w_1$ and $w_2$. 
Thus, $b$ is adjacent in $A_u$ to at least two points $b_1$ and $b_2$. 

We now consider the second statement.
From the definition of the adjacency in $A_u$, the possible segments joining
active points are given in figure \ref{intersecoes}. 
Since the four points are distinct, segments may not meet at endpoints. 
Also, as all points have the same color, the only plausible
intersection of segments are listed in figure \ref{intersecoes}. 
By direct inspection, each such intersection 
violates the definition of adjacency in $A_u$. 
\end{proof}

\begin{proposition}\label{poligonoconvexo} The closure of the 
connected components of $\Phi_b(u)$ are convex polygons with vertices in 
$A_u$ or in the corners of $R$. In particular, $\Phi_b(u)$ is indeed a splitting.
\end{proposition}

\begin{proof} From Lemma \ref{degree}, the juxtaposition of segments in $A_u$, together with
$\partial R$, consists of closed polygonal curves and the vertices of the 
connected components in the statement are either in $A_u$ or in the corners of $R$. 
From the classification of stars,
convexity follows at interior vertices: for boundary vertices, convexity is trivial.
\end{proof}

\subsection{The possible sides and angles of chunks of $\Phi_b(u)$}

Sides of chunks of a splitting $\Phi_b(u)$ are {\emph {active}}  if they consist only of
segments of $A_u$.
{\emph {Inactive sides}} are those which contain no segment of $A_u$.
The remaining sides are called {\emph {mixed}}. Inactive  and mixed sides must lie in 
$\partial R$, since sides in $\Phi_b(u)$ are union of segments in $A_u$ with segments 
in $\partial R$.

\begin{lemma}\label{ladomisto}
There are no mixed sides.
\end{lemma}

\begin{figure}[hbtp]
\centering
\includegraphics[height=0.8in]{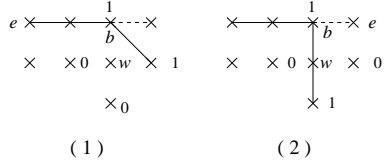}
\caption{Mixed segments do not exist; dashed segments are inactive.}
\label{lado misto}
\end{figure}

\begin{proof}
Suppose $e$ is a mixed side of a chunk in $\Phi_b(u)$. Without loss, suppose $e$ 
belongs to the upper side of $R$. The $m \times n$ rectangular graph $G$
ought to include points immediately below $e$, since $m,n >1$.
Let $b$ be a point of $e$ joining an active and an inactive segment, and let $w$
denote the white point  below $b$.
Figure \ref{lado misto} shows the possible active neighbors of $w$:
any possible choice implies the existence of a segment in $A_u$ meeting $e$
at $b$: in particular, $e$ may not be a side of a chunk.
\end{proof}

\begin{lemma}\label{lado2}
Horizontal or vertical active sides of chunks in $\Phi_b(u)$ have \mbox{length 2}.
\end{lemma}

\begin{proof}
Consider a horizontal active side $h$ of a chunk of $\Phi_b(u)$. 
From the adjacency relation in $A_u$, $h$ has even length. Suppose $h$
has length greater to 2, and take $e$ in $h$ joining two active points at distance 4.

\begin{figure}[hbtp]
\centering
\includegraphics[width=1in]{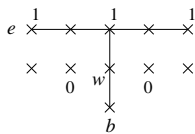}
\caption{Horizontal sides.}
\label{lado que nao and 2}
\end{figure}

As in the previous Lemma, we may suppose that the integral points 
immediately under $h$ belong to $G$. The midpoint of $e$ has at least one
white neighbor $w \in G$, as shown in figure \ref{lado que nao and 2}.
From the adjacency relation in $A_u$, the side neighbors of $w$ are
inactive. There is only one possibility: the other active neighbor of $w$
must be $b$, as in figure \ref{lado que nao and 2}. But then, 
the vertical segment in the same figure ought to be a segment in $A_u$,
contradicting the fact that $e$ lies in a side of a chunk.
\end{proof}

The  angles of chunks are formed by the segments listed in figure \ref{intersecoes} 
combined with segments in $\partial R$. Our next step will be the classification of
angles according to activity of their sides.
Angles defined by an active and an inactive side are {\emph {mixed angles}}. 
From the adjacency relation in $A_u$, angles of  $\pi/4$ in $\Phi_b(u)$ 
are mixed and angles of  $3\pi/4$ are not. We are left with considering right angles.

\begin{lemma}\label{tudo ativo x tudo inativo}
Right angles are not mixed. 
\end{lemma}

\begin{figure}[hbtp]
\centering
\includegraphics[height=0.5in]{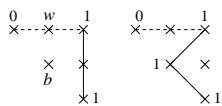}
\caption{An inexistent right angle.}
\label{tudo ativo x tudo inativo fig}
\end{figure}

\begin{proof}
Consider first the case in which the sides of the angle are horizontal and vertical, 
as in figure \ref{tudo ativo x tudo inativo fig}. 
The inactive side must then belong to $\partial R$.
The number of active neighbors of $w$ must be even, and hence $u(b)=1$.
But then, the adjacency relation in $A_u$ implies the presence of 
active segments as in figure \ref{tudo ativo x tudo inativo fig}.

The other kind of right angle has diagonal sides:
these sides may not belong to $\partial R$ and thus have to be active. 
Indeed, by the adjacency relation in $A_u$, this kind of right angle may only 
occur in square chunks of sides with length $\sqrt 2$.
\end{proof}  

\begin{figure}[hbtp]
\centering
\includegraphics[height=0.5in]{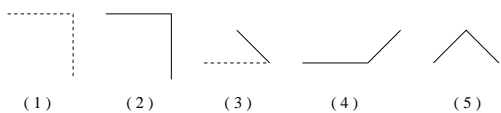}
\caption{Angles; active and inactive sides.}
\label{angulos2}
\end{figure}

Adding up, figure \ref{angulos2} lists all possible angles of chunks in $\Phi_b(u)$, taking
into account the activity of each side.

\begin{lemma}\label{num impar}
The number of $\pi/4$ angles in a chunk of a splitting $\Phi_b(u)$ is even.
\end{lemma}

\begin{proof}
Lemma \ref{ladomisto} excludes the possibility of mixed sides. 
Clearly, every chunk ought to have an even number of mixed angles, and 
the only mixed angles measure $\pi/4$.
\end{proof} 

\begin{lemma}\label{90+90=quadrado}
Let $C$  be a chunk in $\Phi_b(u)$ with a side $s$ common to two right angles 
which is active and either horizontal or vertical.
Then $C$ is a square with sides of length 2.
\end{lemma}

\begin{figure}[hbtp]
\centering
\includegraphics[height=0.8in]{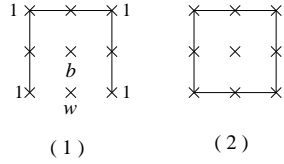}
\caption{(1) implies (2).}
\label{lema 3}
\end{figure}

\begin{proof}
Lemma \ref{lado2} implies that the side $s$ has length 2.
In figure \ref{lema 3}.1, $b$ is not active, again in accordance with the adjacency
relation in $A_u$. Thus, the only active neighbors of $w$ are adjacent in $A_u$.
\end{proof}

\subsection{Chunks in $\Phi_b(u)$ are tiles}

Since the smallest external angle of a chunk is $\pi /4$ 
(figure \ref{angulos2}) chunks of $\Phi_b(u)$ have at most $8$ sides --- we now
classify chunks. \vspace{0.1in}

\noindent{\emph{Octagons:}} 
All inner angles have to measure $3\pi/4$. 
Horizontal and vertical sides have length $2$ since they are active 
(Lemma \ref{lado2}) and thus octagons must be as in figure \ref{pecas}. 
\vspace{0.1in}

\begin{figure}[hbtp]
\centering
\includegraphics[height=0.9in]{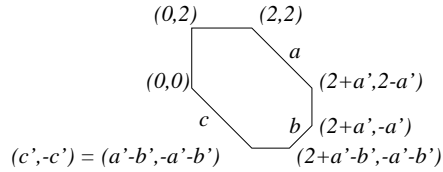}
\caption{There are no heptagonal chunks:
$a=\sqrt2 a'$, $b=\sqrt2 b'$ and $c=\sqrt2 c'.$}
\label{heptagono}
\end{figure}

\noindent{\emph{There are no heptagonal chunks:}} 
A heptagonal chunk should have 6 angles measuring
$3\pi/4$ and one, $\pi/2$ (figure \ref{heptagono}).
Horizontal and vertical sides must be active and have length 2.
From figure \ref{heptagono}, adding coordinates on both axis,  $a'-b'-c'=0$ and
$-a'-b'+c'=0$. Thus, $a-b-c=0$ and $-a-b+c=0$, and hence $b=0$: one side is gone.
\vspace{0.1in}

\noindent{\emph{Hexagons:}} From Lemma \ref{num impar}, chunks of $\Phi_b(u)$
have an even number of angles measuring $\pi/4$.  
Thus, hexagons must have four $3\pi/4$ angles and two $\pi/2$ angles. 
Angles ought to be ordered as in the figure \ref{pecas}: from Lemma \ref{90+90=quadrado}, 
the right angles have to be separated.
\vspace{0.1in}

\noindent{\emph{There are no pentagonal chunks:}} Again, from Lemma \ref{num impar},
 a pentagon should have two $3\pi/4$ angles and three right angles. But then two 
consecutive right angles are inevitable, contradicting
Lemma \ref{90+90=quadrado}. Notice that the hypothesis of the Lemma is
fulfilled: there are no mixed angles and hence no inactive sides.
\vspace{0.1in}

 \noindent{\emph{Quadrilaterals:}} There are two cases.
First consider rectangular chunks. At least one of its sides must be active, otherwise
$u=0$. Since right angles are not mixed, all sides must be active. But then, from
Lemma \ref{lado2} and the adjacency relation in $A_u$, the chunk must be 
one of the squares in figure \ref{pecas}. Second,
the chunk may have two $3\pi/4$ angles and two  $\pi/4$ angles and must be the trapezoidal
tile in figure \ref{pecas}.
\vspace{0.1in}

\noindent{\emph{Triangles:}} The triangle must be right isosceles 
(figure \ref{pecas})
\vspace{0.1in}

Trapezoidal and triangular chunks have dashed sides: 
$\pi/4$ angles are mixed. The upshot of the classification is the following ---
{\emph{the possible chunks in a splitting $\Phi_b(u)$ are the tiles in figure \ref{pecas}:
$\Phi_b(u)$ is a tiling}}.
Thus, the function $\Phi_b$ indeed takes nonzero vectors in the kernel of $BW$ 
to the subset $\mathcal{T}_b$ of tilings. There is an analogous function $\Phi_w$ 
from $\ker \WB \setminus \{0\}$ to $\mathcal{T}_w.$ 

\section{$\Phi_b$ is a bijection}

We start showing that $\Phi_b: \ker BW \setminus \{0\} \to \mathcal{T}_b$ is injective.
Take $u,v \in \ker BW$, inducing tilings $\Phi_b(u)$ and $\Phi_b(v)$.
We have to show that if $\Phi_b(u) = \Phi_b(v)$, then $A_u=A_v$ (and then, trivially, $u = v$).
If $\Phi_b(u)=\Phi_b(v)$, then $A_u$ and $A_v$ may only differ on $\partial R$,
since the only sides which may be active or not lie there. It is
the shape of the tile which  determines which side is active or not: inactive sides are
the large basis of trapezoidal tiles and short sides of right triangles. Thus,
sides of a tile in $\Phi_b(u) = \Phi_b(v)$ belonging to $\partial R$ 
must be simultaneously active or not in both $u$ and $v$, and the proof of injectivity
is complete.

We now consider the surjectivity of $\Phi_b$. 
Given a tiling $T\in \mathcal{T}_b$ (and in particular, the description of type of 
activity for each side), we have to find a nonzero vector $u \in \ker BW$ 
for which $\Phi_b(u)=T$. The sequence of steps below is natural.

\noindent $\bullet$ Defining $u$ from $T$.

\noindent $u(b)=1 \Leftrightarrow b$ belongs to a non-dashed side in $T$.  
\smallskip

\noindent $\bullet$ Checking that $u \in \ker BW$.

\noindent We must show that, for each white point $w$, 
the sum of the values $u(b)$ for its neighbors $b \in G$ is even. By checking
for each kind of tile, it is clear that the sum of $u(b)$ for neighbors $b$ of
$w$ in a single tile is even. Thus, if $w$ lies in the interior of a tile, we
are done. Suppose now that $w$ lies in a side of a tile of $T$.
First, notice that vertices of tiles belonging to active sides are black. Hence,
if $w$ is such a vertex, it must belong to inactive sides, and thus, both sides
containing $w$ lie in $\partial R$ --- $w$ must be a corner of $R$,
belonging to a right triangle. In this case, the black neighbors of $w$ are active,
if the short sides of the triangle are of length 1, or inactive, for longer lengths.
In both cases, the total sum of $u(b)$ over black neighbors $b$ of $w$ is even.

Finally, if $w$ belongs to a side of a tile without being a vertex, this side must
be either horizontal or vertical, since diagonal sides only cross black points.
If this side lies in $\partial R$, all the black neighbors of $w$ 
belong to the same tile, and we are done. If, instead, this side is in the interior
of $R$, we may take it to be horizontal and it must be an active side.
Thus, $u(b_1)=u(b_2)=1$, where $b_1$ and $b_2$ are the side neighbors of $w$.
Besides, the black neighbors of $w$ belong to two tiles, since $w$ is not a vertex.
Suppose $b_1$, $b_2$ and $b_3$ in one tile and $b_1$, $b_2$ and $b_4$ in another.
Then
$u(b_1)+u(b_2)+u(b_3)$ and 
$u(b_1)+u(b_2)+u(b_4)$ are even.
Since $u(b_1)=u(b_2)=1$, we must have \mbox{$u(b_3)=u(b_4)=0$}. 
Hence, $u(b_1)+u(b_2)+u(b_3)+u(b_4)=0 \pmod{2}.$
\smallskip

\noindent $\bullet$ The vector $u$ obtained from $T$ indeed satisfies $\Phi_b(u)=T$. 

\noindent The black points by which non-dashed sides of $T$ and $\Phi_b(u)$ passes are the same.
Indeed, these are the active points of $u$. We are left with checking that 
the segments between these points are the same in $T$ and $\Phi_b(u)$. 
Clearly, $\partial R$ consists of segments common to both tilings. Therefore dashed 
sides are irrelevant: we only need to prove that $T$ and $\Phi_b(u)$ have the same active
segments.

\begin{figure}[hbtp]
$$
\begin{array}{ccccc}
\includegraphics[width=0.7in]{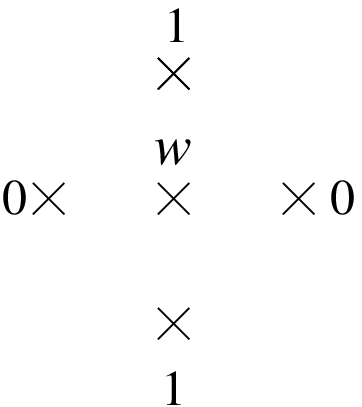}&\hspace*{0.5cm}&
\includegraphics[width=0.8in]{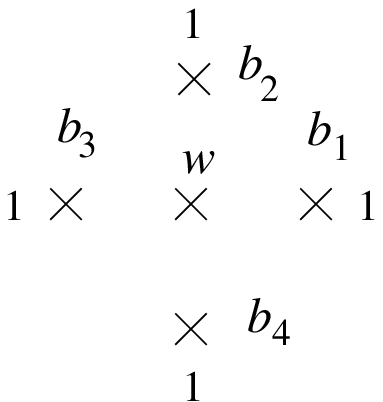}&\hspace*{0.5cm}&
\includegraphics[width=0.7in]{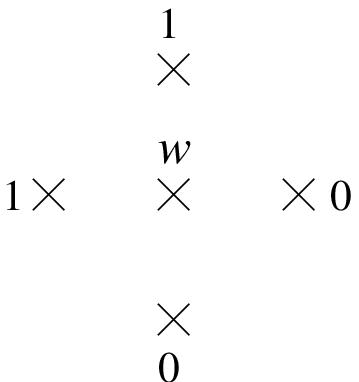}\\\\
( 1 )&&( 2 )&&( 3 )\\\\
\includegraphics[width=0.7in]{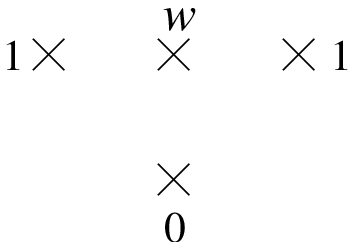}&\hspace*{0.5cm}&
\includegraphics[width=0.4in]{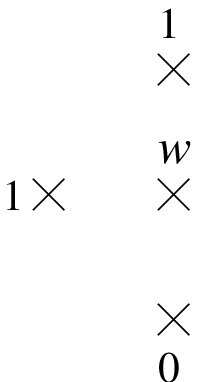}&\hspace*{0.5cm}&
\includegraphics[width=0.4in]{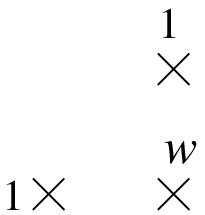}\\\\
( 4 )&&( 5 )&&( 6 )
\end{array}
$$
\caption{Nontrivial active neighborhoods.}
\label{vizinhancas}
\end{figure}

The {\emph{active neighborhood}} of a point $w \in G$ under $u \in \ker BW$ 
is the set of its active black neighbors in $G$. An active neighborhood is
{\emph{trivial}} if there are no active points. A simple check obtains all nontrivial active 
neighborhoods: they are listed in figure \ref{vizinhancas}, up to rotations. 
In cases (1), (2) and (3), $w$ is an interior point; in (4) and (5), 
$w$ belongs to a side of $\partial R$, and in (6), $w$ is a corner of $R$. 
Trivial active neighborhoods are not relevant in the forthcoming argument.
In order to show that $T$ and $T_\Phi$ have the same active segments, all we have to do is 
to show this fact for each active neighborhood.

\begin{figure}[h]
\centering
\includegraphics[height=0.4in]{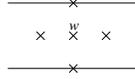}
\caption{The boundary of tile $t$.}
\label{situacao1}
\end{figure}

In the case of a neighborhood of type (1), the two active points are joined by a segment in
$\Phi_b(u)$, by definition. Suppose that this segment does not belong to $T$: 
$w$ then is in the interior of a tile $t$ of $T$. The two active neighbors of
$w$ belong to the boundary of $t$, by the construction of $u$ from $T$.
Since tiles are convex, the segments of $T$ form angles which are less than or equal to 
$\pi$. Also, sides of tiles in $T$ do not pass by points at which $u$ is zero 
(we do not have to take into account boundary points). 
Thus, the sides of $t$ passing by the active points of the neighborhood of $w$ must
be horizontal, as shown in figure \ref{situacao1}; these sides extend to active neighbors
at distance 2. But, by figure \ref{pecas} horizontal sides have length at most 2 
-- contradiction.
Thus $T$ equals $\Phi_b(u)$, when restricted to an active neighborhood of type (1).

\begin{figure}[hbtp]
\centering
\includegraphics[height=0.7in]{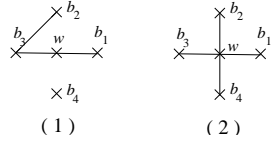}
\caption{Active segments.}
\label{situacoes3}
\end{figure}

For a neighborhood of type (2), $T$ contains active segments joining each neighbor of $w$ to some
other neighbor. Consider the possibility of one of this segments being horizontal,
as $\overline{b_1b_3}$ in figure \ref{situacoes3}.1. Neighbor $b_2$ is either joined to $b_3$ 
or to $b_4$, up to trivial symmetries, as in figure \ref{situacoes3}.
Tiling $T$ cannot contain (1), since the
$\pi/4$ angle may only appear at $\partial R$ (again, check the shape of tiles
in figure \ref{pecas}).
Possibility (2) may not happen either, since it implies a white vertex of a tile which is 
not a corner of $R$. Thus segment $\overline{b_1b_3}$ is not in $T$. Vertical segments
are equally excluded. 

Suppose now that there are only diagonal active segments. The only tile with
parallel diagonal sides $\sqrt2$ apart is the tilted square. Thus, the four diagonal
segments must appear in both $T$ and $\Phi_b(u)$.

The study of neighborhoods of types (3), (4), (5) and (6) is similar and will be
left to the reader. This proves the surjectivity of $\Phi_b$.

At this point, the proof of Theorem  \ref{principal} is complete.

\section{Counting polarized harmonic functions}

To prove Theorem \ref{formula}, we compute $\# \mathcal{T}$, the number of tilings of $R$. 
Clearly,
\begin{align*}
\# \mathcal{T} &= \# \mathcal{T}_b \cup \mathcal{T}_w 
&&\text{(by Proposition \ref{coordenadas inteiras})} \\
    &= \#\mathcal{T}_b \, + \,  \# \mathcal{T}_w   
&&\text{(since $\mathcal{T}_b$ and $\mathcal{T}_w$ are disjoint)} \\
    &= \# \ker BW\setminus \{0\} \, + \, \# \ker \WB \setminus\{0\}
                       &&\text{(by Theorem \ref{principal})}      \\
    &= \#\ker BW  \, + \, \# \ker \WB -2.
\end{align*}

\begin{figure}[hbtp]
\centering
\includegraphics[height=1.8in]{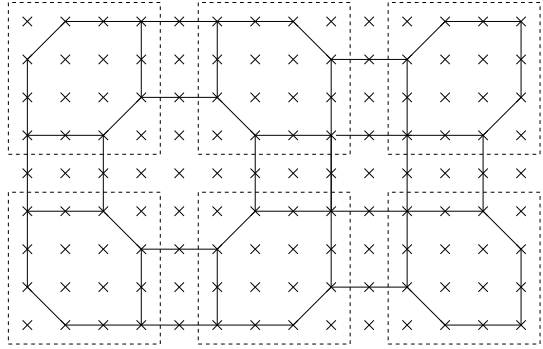}
\caption{Grid, fundamental squares and symmetries.}
\label{espelhamentos}
\end{figure}

We only have to compute the dimensions $\beta$ and $\omega$ of the kernels of $BW$ and $\WB$.
We will exhibt two symmetries common to the kernel elements, and 
which are evident in figure \ref{espelhamentos}. 
Given a vector $v$ belonging to $(\ZZ/(2))^{mn}$, denote by $\sigma_{i,j}$ the sum $\pmod 2$ 
of the coordinates of $v$ at the four neighbors of $(i,j)$ in $G$:
\mbox{$\sigma_{i,j}=v_{i,j-1}+v_{i,j+1}+v_{i-1,j}+v_{i+1,j} \pmod 2$}. 
Here  $v_{i,j}=0$ if $(i,j)\not\in G.$
Recall that the kernel of $M$, the adjacency matrix of $G$, is the direct sum of the
kernels of $BW$ and $\WB$.
Without loss we will 
assume from now on that the rectangular graph satisfies $m \geq n$.

\begin{proposition}\label{simetria}
Let $v \in \ker M$. Then 
\begin{itemize}
\item [(i)] $v_{i,j}=v_{j,i},$ if both points belong to $G$.
\item [(ii)] If the vector $v$ equals zero along the column $x=m_0$
then
\mbox{$v_{m_0 - i,j} = v_{m_0+i,j},$} again if both points belong to $G$.
\end{itemize}
\end{proposition}

Clearly, there are analogous symmetries across the other diagonals of the graph, 
and across rows of zeros of $v$.

\begin{proof}
To start off the induction argument, notice that $v_{i,j} =v_{j,i}$, if $i = j$ (i.e., 
along the line $y = x$). Now compare values along the lines
$y = x+1$ and $y=x-1$. Since $\sigma_{0,0} = 0,$ we must have
$v_{1,0} = v_{0,1}$. Sequentially,
$\sigma_{1,1}=0,\ \sigma_{2,2}=0,\ldots$ imply in turn
$v_{2,1} = v_{1,2}$, \mbox{$\ v_{3,2} = v_{2,3},\ldots$} 
Compare now values along the lines $y = x+2$ and $y = x-2$, making use of
$\sigma_{1,0} = 0 = \sigma_{0,1}$. Proceed until the lines $y = x \pm k$
cover the square $[0,n-1] \times [0,n-1]$.

The proof of (ii) is a simple induction on $i$. Case $i=1$ is as hard as
the general case, and follows from expanding
$0 = \sigma_{m_0,j} = v_{m_0-1,j} + v_{m_0+1,j}.$
\end{proof}

Consider the $m \times n$ rectangular graph $G$  and let
$c = \gcd(m+1,n+1)-1$. The {\emph{grid}} $\mathcal{G}$ of $G$ is the set of points in
$G$ belonging to the union of the horizontal lines 
$y = k(c+1)-1$ with the vertical lines $x = k(c+1)-1$,  for integer $k$.
The graph $G$, after removal of its grid $\mathcal{G}$, splits into $c \times c$ squares,
called {\emph{fundamental squares}}.
Figure \ref{espelhamentos} shows the grid and the fundamental squares for the 
$14\times 9$ rectangular graph.

\begin{theorem}\label{engradado} 
Let $G$ be the $m \times n$ rectangular graph with grid $\mathcal{G}$ 
and adjacency matrices $BW_{m,n}$ and $\WB_{m,n}$, with kernel dimensions
$\beta$ and $\omega$ respectively.
Let $v \in \ker BW_{m,n}$ (resp. $\ker \WB_{m,n}$). Then
\begin{itemize}
\item[(i)] The vector $v$ is identically zero on $\mathcal{G}$.
\item[(ii)]
The restriction of $v$  to each fundamental square lies in the kernel of the matrix
$BW_{c,c}$ (resp. $\WB_{c,c}$) of this square.
\item[(iii)] Let $s_1$ and $s_2$ be two adjacent fundamental squares, bounded by a common 
row or column of $\mathcal{G}$. Then the restrictions of $v$ to $s_1$ and $s_2$ are 
mirror images of each other with respect to this common row or column.
\item[(iv)] $\beta = \dim \ker BW_{c,c} \text{ and } 
\omega = \dim \ker \WB_{c,c}.$
\end{itemize}
\end{theorem}

\begin{proof}
We shall prove these assertions by induction on the set of pairs $(m,n)$ for which 
\mbox{$\gcd(m+1,n+1) = c + 1$} is fixed, ordered lexicographically.
For the starting point $m = n =c$, the result is obvious.
Suppose that \mbox{$v^{m,n} \in \ker BW_{m,n}$}, in a self-explanatory notation.
From the symmetry (i) of the previous proposition, the restriction of $v^{m,n}$ to 
column $x = n$ is zero. Indeed, for $v = v^{m,n}$, since
\begin{alignat*}{3}
0 &= \sigma_{n-1,j} &&= v_{n-2,j} + v_{n-1,j-1} + v_{n-1,j+1},\\
0 &= \sigma_{j,n-1} &&= v_{j,n-2} + v_{j-1,n-1} + v_{j+1,n-1} + v_{j,n},
\end{alignat*}
we must have $v_{j,n} = 0$. Removal of column $x = n$  from $G = G_{m,n}$ gives rise
to two connected subgraphs: a square $G_{n,n}$ containing the first $n$ columns of $G$ and 
a possibly rectangular subgraph $G_{m-n-1,n}$ containing only the last $m-n-1$ columns. 
Also, by making use of the column of zeros and the fact that $v^{m,n} \in \ker BW_{m,n}$, 
we have that the restriction
$v^{m-n-1,n}$ of $v^{m,n}$ to $G_{m-n-1,n}$ belongs to $\ker BW_{m-n-1,n}$. 
Notice that \mbox{$\gcd(m-n-1+1,n+1) = c+1$}: by induction, then, $v^{m-n-1,n}$ 
has zeros on the grid of $G_{m-n-1,n}$. Thus, the nonzero entries of
$v^{m-n-1,n}$ belong to $c \times c$ squares framed by points in which 
$v^{m-n-1,n}$ is zero. Now using symmetry (ii), mirroring these squares across columns
$n$, $n-c-1$ \ldots, we see that the square  $G_{n,n}$ also splits in 
$c \times c$ squares with the same property. Conversely, mirroring a 
kernel element on the fundamental square across the grid obtains a kernel
element on the full rectangle: this proves (iv).
\end{proof}

We now compute the dimensions of the kernels of $BW$ and $\WB$ for the square graph $G_{c,c}$.

\begin{proposition}\label{quadrado cheio}                                                      
Let $G_{c,c}$ be the $c \times c$ square graph. Then 
$$2\dim \ker BW_{c,c} = c + (c \bmod 2)\quad \text{ and } \quad
2\dim \ker \WB_{c,c} = c - (c \bmod 2).$$

\end{proposition}

\begin{proof}
It is clear that a vector in $\ker BW_{c,c}$ (or in $\ker \WB_{c,c}$) 
is determined by its values along the first column $x = 0$ of $G_{c,c}$. The fact
that any choice of values along this columns indeed gives rise to a kernel vector 
follows from the diagonal symmetry (i) of Proposition \ref{simetria}.
\end{proof}

\begin{figure}[hbtp]          
\centering
\includegraphics[height=2in]{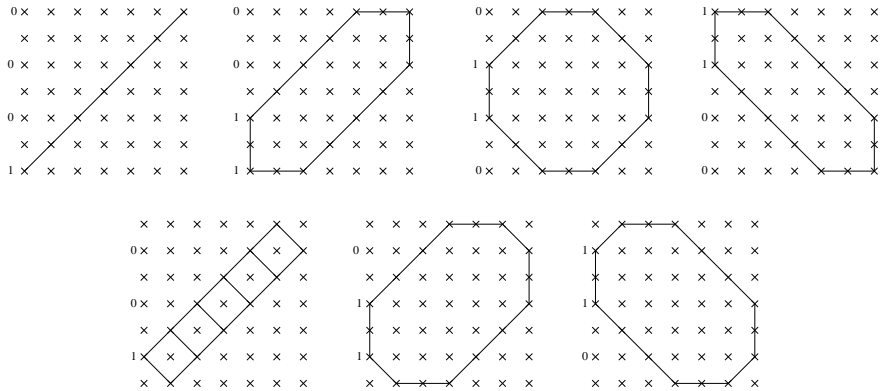}
\caption{Basis for $\ker BW$ and for $\ker \WB$.}
\label{nucleo impar}
\end{figure}

Certain choices of values along the first column of $G_{c,c}$ obtain 
especially simple kernel vectors. In figure \ref{nucleo impar},
the first (resp. second) row indicates a basis for $\ker BW$ (resp. $\ker \WB$).

\begin{figure}[h]
\centering
\includegraphics[height=5.6in]{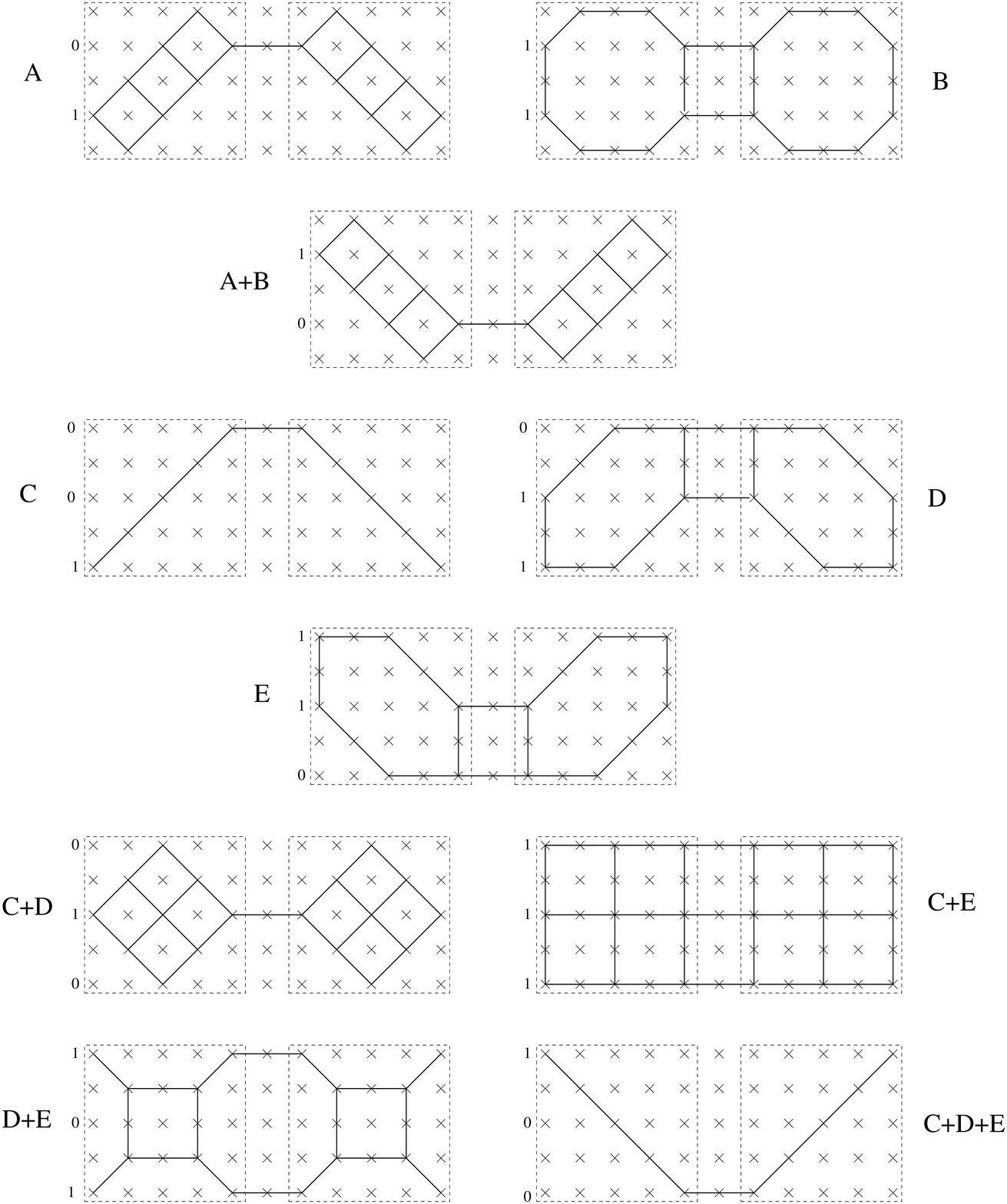}
\caption{All tilings of the $10 \times 4$ rectangle.}
\label{metodo}
\end{figure}

Combining the previous proposition with item (iv) of the theorem above,
we complete the proof of Theorem  \ref{formula}.
Figure \ref{metodo} shows all tilings of the $10 \times 4$ rectangle:
the nonzero elements in the kernel of $\WB$ give rise to the first three, $A$, $B$ and $A+B$.
The remaining ones correspond to nonzero elements in the kernel of $BW$.

%
%
%

\end{document}